\documentclass[12pt]{amsart}
\usepackage{amsfonts, amssymb, mathrsfs, amsmath}
\usepackage{fullpage, verbatim, bbm}
\usepackage[colorlinks=true]{hyperref}
\usepackage{breakurl}
\usepackage{tikz}
\usepackage{multirow,tabularx}
\usepackage{xfrac}
\usepackage{longtable}
\usepackage[bb=libus]{mathalpha}

\usepackage{etoolbox}
\makeatletter
\let\ams@starttoc\@starttoc
\makeatother
\usepackage[parfill]{parskip}
\makeatletter
\let\@starttoc\ams@starttoc
\patchcmd{\@starttoc}{\makeatletter}{\makeatletter\parskip\z@}{}{}
\makeatother

\newcommand\CA{{\mathscr A}} 
\newcommand\Ac{{\mathscr A}} 
\newcommand\CB{{\mathscr B}}
 
\newcommand\cC{{\mathcal C}}

\newcommand\CAF{{\mathcal {AF}}} 
\newcommand\CIF{{\mathcal {IF}}} 
 
\newcommand\CMF{{\mathcal {MF}}}

\newcommand\BBk{{\mathbb k}}

\newcommand\BBZ{{\mathbb Z}}

\newcommand\BBF{{\mathbb F}}

\newcommand\codim{\operatorname{codim}}

\newcommand\Der{{\operatorname{Der}}}

\newcommand\rank{\operatorname{rank}}
\newcommand\rk{\operatorname{rk}}

\numberwithin{equation}{section}

\theoremstyle{plain}

\newtheorem{theorem}[equation]{Theorem}

\theoremstyle{definition}
\newtheorem{defn}[equation]{Definition}
\newtheorem{remark}[equation]{Remark}
\newtheorem{remarks}[equation]{Remarks}
\newtheorem{example}[equation]{Example}

\begin{document}

\title[MAT-Freeness is not combinatorial]
{MAT-Freeness  is not combinatorial}

\author[T.~Hoge and G.~R\"ohrle]{Torsten Hoge and Gerhard R\"ohrle}
\address
{Fakult\"at f\"ur Mathematik,
	Ruhr-Universit\"at Bochum,
	D-44780 Bochum, Germany}
\email{torsten.hoge@rub.de}
\email{gerhard.roehrle@rub.de}


\allowdisplaybreaks

\begin{abstract}
We show that the notion of
MAT-freeness for hyperplane arrangements depends on the underlying field.
In particular, MAT-freeness is not combinatorial.
\end{abstract}

\maketitle



\section{Introduction}
 
In the theory of hyperplane arrangements the study of classes which are determined purely by combinatorial data, i.e., which only depend on their underlying intersection lattices (see Definition \ref{def:comb}),  has a  long tradition.

Undoubtedly, the most tantalizing example  is the question whether freeness for hyperplane arrangements is combinatorial or not. The assertion that this is the case is known as \emph{Terao's conjecture} \cite[p.~293]{terao:freeI}. 
In 1990 Ziegler showed that the question of freeness does depend on the choice of the underlying field, \cite[Ex.~4.1]{ziegler}. Thus freeness is not combinatorial.
Subsequently, in the literature Terao's conjecture specifies a fixed underlying field, \cite[Conj.~4.138]{orlikterao:arrangements}.
While overwhelming evidence in support of this conjecture stems from Terao's seminal Factorization Theorem for the Poincar\'e polynomial of a 
free arrangement \cite[Thm.~4.137]{orlikterao:arrangements}, 
a conclusive resolution eludes until this day. 

There are stronger notions of freeness, such as \emph{inductive freeness} (Definition \ref{def:indfree-simple}) or \emph{additive freeness} (Definition \ref{def:addfree}). In contrast to freeness itself, both of  them are combinatorial, see \cite[Lem.~2.5]{cuntzhoge}, respectively  \cite[Thm.~1.4]{abe:sf}.

In  \cite{cuntzhoge} Cuntz and Hoge showed that the weaker notion of \emph{recursive freeness} (\cite[Def.~4.60]{orlikterao:arrangements}) is again not combinatorial. 
In this short note we observe that the related concept of
\emph{MAT-freeness} (Definition \ref{def:mf}), due to Cuntz and  M\"ucksch \cite{CunMue19_MATfree}, also fails to be combinatorial: Example \ref{ex:nonmf} shows that MAT-freeness depends on the choice of the underlying field. Following \cite{CunMue19_MATfree}, we denote the class of MAT-free arrangements by $\CMF$.

\begin{theorem}
	\label{thm:main}
	\begin{itemize}
		\item [(i)] The class $\CMF$ is not combinatorial.
		\item [(ii)] The subclass of $\CMF$ consisting of arrangements over infinite fields is combinatorial. 
	\end{itemize}
\end{theorem}

The origin of MAT-freeness stems from the so called \emph{ Multiple Addition Theorem} \ref{Thm_MAT} from \cite{ABCHT16_FreeIdealWeyl}. 
In \emph{loc.~cit.}  and \cite{CunMue19_MATfree} the investigations concerned hyperplane arrangements  over fields of characteristic $0$. MAT-freeness is indeed combinatorial for the subclass consisting of hyperplane arrangements defined over such fields, 
as shown in \cite[Lem.~3.4]{CunMue19_MATfree}; 
see also Remark \ref{rem:main}.

For general information on hyperplane arrangements we refer to  \cite{orlikterao:arrangements}.

\section{Preliminaries}
\label{sect:prelim}

\subsection{Hyperplane arrangements}
\label{ssect:hyper}

Let $\BBk$ be a field and let $V = \BBk^\ell$ 
be an $\ell$-dimensional $\BBk$-vector space.
A \emph{hyperplane arrangement} is a pair
$(\CA, V)$, where $\CA$ is a finite collection of hyperplanes in $V$.
Usually, we simply write $\CA$ in place of $(\CA, V)$.
We write $|\CA|$ for the number of hyperplanes in $\CA$.
The empty arrangement in $V$ is denoted by $\Phi_\ell$.

The \emph{lattice} $L(\CA)$ of $\CA$ is the set of subspaces of $V$ of
the form $H_1\cap \dotsm \cap H_i$ where $\{ H_1, \ldots, H_i\}$ is a subset
of $\CA$. 
For $X \in L(\CA)$, we have two associated arrangements, 
firstly
$\CA_X :=\{H \in \CA \mid X \subseteq H\} \subseteq \CA$,
the \emph{localization of $\CA$ at $X$}, 
and secondly, 
the \emph{restriction of $\CA$ to $X$}, $(\CA^X,X)$, where 
$\CA^X := \{ X \cap H \mid H \in \CA \setminus \CA_X\}$.
The lattice $L(\CA)$ is a partially ordered set by reverse inclusion:
$X \le Y$ provided $Y \subseteq X$ for $X,Y \in L(\CA)$.
The rank function on $L(\CA)$ is given by $\rank(X)=\codim X$ for $X\in L(\CA)$.
The rank of $\CA$, $\rank(\CA)$,  is the rank of a maximal element of $L(\CA)$.

\begin{defn}
	\label{def:comb}
Let $\cC$ be a class of arrangements and let $\CA \in \cC$. 
		If for all arrangements $\CB$ with $L(\CB)$ and $L(\CA)$ isomorphic as lattices, we have $\CB \in  \cC$, then the class $\cC$ is called \emph{combinatorial}.
		Here $\CA$ and $\CB$ do not have to be defined over the same underlying field.
\end{defn}

\subsection{Free arrangements}
\label{ssect:freeness}
Let $S = S(V^*)$ be the symmetric algebra of the dual space $V^*$ of $V$.
If $x_1, \ldots , x_\ell$ is a basis of $V^*$, then we identify $S$ with 
the polynomial ring $\BBk[x_1, \ldots , x_\ell]$.
Letting $S_p$ denote the $\BBk$-subspace of $S$
consisting of the homogeneous polynomials of degree $p$ (along with $0$),
$S$ is naturally $\BBZ$-graded: $S = \oplus_{p \in \BBZ}S_p$, where
$S_p = 0$ in case $p < 0$.

Let $\Der(S)$ be the $S$-module of algebraic $\BBk$-derivations of $S$.
Using the $\BBZ$-grading on $S$, $\Der(S)$ becomes a graded $S$-module.
For $i = 1, \ldots, \ell$, 
let $D_i := \partial/\partial x_i$ be the usual derivation of $S$.
Then $D_1, \ldots, D_\ell$ is an $S$-basis of $\Der(S)$.
We say that $\theta \in \Der(S)$ is 
\emph{homogeneous of polynomial degree p}
provided 
$\theta = \sum_{i=1}^\ell f_i D_i$, 
where $f_i$ is either $0$ or homogeneous of degree $p$
for each $1 \le i \le \ell$.
In this case we write $\deg \theta = p$.

Let $\CA$ be an arrangement in $V$. 
Then for $H \in \CA$ we fix $\alpha_H \in V^*$ with
$H = \ker(\alpha_H)$.
The \emph{defining polynomial} $Q(\CA)$ of $\CA$ is given by 
$Q(\CA) := \prod_{H \in \CA} \alpha_H \in S$.

The \emph{module of $\CA$-derivations} of $\CA$ is 
defined by 
\[
D(\CA) := \{\theta \in \Der(S) \mid \theta(\alpha_H) \in \alpha_H S
\text{ for each } H \in \CA \} .
\]
We say that $\CA$ is \emph{free} if the module of $\CA$-derivations
$D(\CA)$ is a free $S$-module.

With the $\BBZ$-grading of $\Der(S)$, 
also $D(\CA)$ 
becomes a graded $S$-module,
\cite[Prop.~4.10]{orlikterao:arrangements}.
If $\CA$ is a free arrangement, then the $S$-module 
$D(\CA)$ admits a basis of $\ell$ homogeneous derivations, 
say $\theta_1, \ldots, \theta_\ell$, \cite[Prop.~4.18]{orlikterao:arrangements}.
While the $\theta_i$'s are not unique, their polynomial 
degrees $\deg \theta_i$ 
are unique (up to ordering). This multiset is the set of 
\emph{exponents} of the free arrangement $\CA$
and is denoted by $\exp \CA$.
We use the shorthand notation  
$(e_1,\ldots,e_\ell)_\leq$ for  
$(e_1 \leq e_2 \leq \ldots \leq e_\ell)$ when we address the set of 
exponents of a free arrangement.

The fundamental \emph{Addition Deletion Theorem} 
due to Terao  \cite{terao:freeI} plays a 
crucial role in the study of free arrangements, 
\cite[Thm.~4.51]{orlikterao:arrangements}. 
Suppose $\CA \neq \Phi_\ell$. Fix a member $H_0$ in $\CA$. Set $\CA' = \CA \setminus \{H_0\}$ and  $\CA'' = \CA^{H_0}$. Then $(\CA, \CA', \CA'')$ is frequently referred to as a \emph{triple of arrangements} (with respect to $H_0$). 

\begin{theorem}
	\label{thm:add-del-simple}
	Suppose $\CA \neq \Phi_\ell$ and
	let $(\CA, \CA', \CA'')$ be a triple of arrangements. Then any 
	two of the following statements imply the third:
	\begin{itemize}
		\item[(i)] $\CA$ is free with $\exp\CA = (b_1, \ldots , b_{\ell -1}, b_\ell)$;
		\item[(ii)] $\CA'$ is free with $\exp\CA' = (b_1, \ldots , b_{\ell -1}, b_\ell-1)$;
		\item[(iii)] $\CA''$ is free with $\exp\CA'' = (b_1, \ldots , b_{\ell -1})$.
	\end{itemize}
\end{theorem}

Theorem \ref{thm:add-del-simple} motivates the following notion.

\begin{defn}
	[{\cite[Def.~4.53]{orlikterao:arrangements}}]
	\label{def:indfree-simple}
	The class $\CIF$ of \emph{inductively free} arrangements 
	is the smallest class of arrangements subject to
	\begin{itemize}
		\item[(i)] $\Phi_\ell \in \CIF$ for each $\ell \ge 0$;
		\item[(ii)] if there exists a hyperplane $H_0 \in \CA$ such that both
		$\CA'$ and $\CA''$ belong to $\CIF$, and $\exp \CA '' \subseteq \exp \CA'$, 
		then $\CA$ also belongs to $\CIF$.
	\end{itemize}
\end{defn}

\begin{remark}
	\label{rem:IF}
	The class $\CIF$ is known to be combinatorial \cite[Lem.~2.5]{cuntzhoge}.
\end{remark}

\subsection{Additively free arrangements}

In view of Theorem \ref{thm:add-del-simple}, it is also natural to consider the following property, \cite[Def.~1.6]{abe:sf} which is weaker than inductive freeness. 

\begin{defn}
	\label{def:addfree}
	An arrangement $\CA$ is called 
	\emph{additively free} if there is a \emph{free filtration} 
	\[
	\Phi_\ell = \CA_0 \subset \CA_1 \subset \cdots \subset \CA_n = \CA,
	\]
	of $\CA$, i.e., where each $\CA_i$ is free with $|\CA_i| = i$ for each $i$. 
	Denote this class by $\CAF$.
\end{defn}

\begin{remarks}
	\label{rem:AF}
	
	(i).
	In view of  \cite[Thm.~4.46]{orlikterao:arrangements}, we have $\CIF \subseteq \CAF$.
	However, there are examples of  arrangements  which are additively free but not inductively free, see \cite{hogeroehrle:stairfree}.
	
	(ii). The class $\CAF$ is known to be combinatorial \cite[Thm.~1.4]{abe:sf}.
\end{remarks}

\subsection{MAT-free arrangements}
\label{SSec_MAT}

Next we consider a freeness concept due to Cuntz and M\"ucksch, \cite{CunMue19_MATfree} which is closely related to inductive freeness. 
We begin by recalling the core result from \cite{ABCHT16_FreeIdealWeyl}, the so called \emph{Multiple Addition Theorem}.
This can be viewed as a generalization of the addition part of  Theorem
\ref{thm:add-del-simple}.

\begin{theorem}[{\cite[Thm.~3.1]{ABCHT16_FreeIdealWeyl}}] 
\label{Thm_MAT}
Let $\Ac' = (\Ac', V)$ be a free arrangement with
$\exp(\Ac')=(e_1,\ldots,e_\ell)_\le$
and let $1 \le p \le \ell$ be the multiplicity of the highest exponent, i.e.
\[ e_{\ell-p} < e_{\ell-p+1} =\cdots=e_\ell=:e. \]
Let $H_1,\ldots,H_q$ be hyperplanes in $V$ with
$H_i \not \in \Ac'$ for $i=1,\ldots,q$. Define
\[ \Ac''_j:=(\Ac'\cup \{H_j\})^{H_j}=\{H\cap H_{j} \mid H\in \Ac'\}, \quad \text{ for }j=1,\ldots,q. \]
Assume that the following conditions are satisfied:
\begin{itemize}
	\item[(1)]
	$X:=H_1 \cap \cdots \cap H_q$ is $q$-codimensional.
	\item[(2)]
	$X \not \subseteq \bigcup_{H \in \Ac'} H$.
	\item[(3)]
	$|\Ac'|-|\Ac''_j|=e$ for $1 \le j \le q$.
\end{itemize}
Then $q \leq p$ and $\Ac:=\Ac' \cup \{H_1,\ldots,H_q\}$ is free with
$$\exp(\Ac)=(e_1,\ldots,e_{\ell-q},e+1,\ldots,e+1)_\le.$$
\end{theorem}

We consider the simultaneous addition of several hyperplanes using
Theorem \ref{Thm_MAT}. This motivates the next terminology.

\begin{defn}
Let $\Ac'$, $\{H_1,\ldots,H_q\}$ be as in Theorem \ref{Thm_MAT} such that
conditions (1)--(3) are satisfied. Then the addition 
of $\{H_1,\ldots,H_q\}$ resulting in $\Ac = \Ac' \cup \{H_1,\ldots,H_q\}$
is called an \emph{MAT-step}.
\end{defn}

An iterative application of Theorem \ref{Thm_MAT} motivates the following natural concept.

\begin{defn}%
	[{\cite[Def.~3.2, Lem.~3.8]{CunMue19_MATfree}}]
	\label{def:mf}%
	
	An arrangement $\Ac$ is called \emph{MAT-free} if there exists an ordered partition
	\[
	\pi = (\pi_1|\cdots|\pi_n)
	\] 
	of $\Ac$ such that for every $0 \leq k \leq n-1$ the following hold. Set $$\Ac_k := \bigcup_{i=1}^k \pi_i.$$ Then 
	\begin{itemize}
		\item[(a)] $\rk(\pi_{k+1}) = \vert \pi_{k+1} \vert$,
		\item[(b)] $\cap_{H \in \pi_{k+1}} H  
		\nsubseteq \bigcup_{H' \in \Ac_k}H'$,
		\item[(c)] $\vert \Ac_k \vert - \vert (\Ac_k \cup \{H\})^H \vert = k$ for each $H \in \pi_{k+1}$,
	\end{itemize}
	i.e.\ $\Ac_{k+1} = \Ac_{k}\cup\pi_{k+1}$ is an MAT-step.
	An ordered partition $\pi$ with these properties is called an \emph{MAT-partition} for $\Ac$.
\end{defn}

Following \cite{CunMue19_MATfree}, the class of MAT-free arrangements is denoted by $\CMF$. 

Note that in \cite{CunMue19_MATfree} MAT-freeness was defined differently.
However, for our purpose its characterization in \cite[Lem.~3.8]{CunMue19_MATfree} is sufficient.
Hence we take the latter here for our definition.

\begin{remark}
	\label{rem:MAT-free}
	Suppose that $\Ac$ is MAT-free with MAT-partition
	$\pi = (\pi_1|\cdots|\pi_n)$. Then 
	\begin{enumerate}
		\item [(a)]
		for each $1 \leq k \le n$, $\Ac_k$ is MAT-free with MAT-partition $(\pi_1|\cdots|\pi_{k})$, 
		
		\item [(b)]
		$\Ac$ is free and the exponents 
		$\exp(\Ac) = (e_1,\ldots,e_\ell)_\le$ of $\Ac$ are given by the block sizes
		of the dual partition of $\pi$: 
		\[ 
		e_i := |\{k \mid |\pi_k|\geq \ell-i+1 \}|, 
		\]
		
		\item [(c)]
		$|\pi_1| > |\pi_2| \geq \cdots \geq |\pi_n|$.
	\end{enumerate}
	Statement (a) is clear by definition.
	Statements (b) and (c) follow readily from Theorem \ref{Thm_MAT} and a simple induction.
\end{remark}

\section{Proof of Theorem \ref{thm:main}}

The following example proves Theorem \ref{thm:main}(i).

\begin{example}
	\label{ex:nonmf}
Consider the arrangement $\CA$ over $\BBF_2$ given by the defining polynomial
$$
    Q(\CA) = xyz (x+y) (x+z) (y+z).
$$
Furthermore consider the arrangement $\CB$ over $\BBF_4$ given by the same defining polynomial
$$
    Q(\CB) = xyz (x+y) (x+z) (y+z).
$$
Then
\begin{itemize}
\item[(1)] both arrangements are free with $\exp(\CA) = \exp(\CB) = (1,2,3)$; 
\item[(2)] the lattices $L(\CA)$ and $L(\CB)$ are isomorphic;
\item[(3)] the arrangement $\CB$ is MAT-free but $\CA$ is not MAT-free.
\end{itemize}

Since the coefficients of the linear forms in the defining polynomial $Q(\CA) = Q(\CB)$ are elements of $\BBF_2$, both lattices are isomorphic.
Both arrangements are additively free and in dimension $3$ they are automatically inductively free.
The order of the hyperplanes, as given in the definition of the defining polynomials above yield an inductive sequence of hyperplanes yielding successive inductively free subarrangements for both, $\CA$ and $\CB$.

Let  $\pi =(\{\ker x, \ker y, \ker z\}, \{ \ker (x+y), \ker (x+z)\}, \{\ker (y+z)\})$. One checks that this is an MAT-partition for $\CB$ and hence $\CB$ is MAT-free over $\BBF_4$.

In contrast, $\CA$ is not MAT-free. We argue by contradiction. Suppose that $\CA$ is MAT-free. Then there is a partition $\pi = (\pi_1,\pi_2,\pi_3)$ of $\CA$ such that
$\pi_3$ is a singleton. After possibly renaming the coordinates,  we may assume that $\pi_3$ consists of either $\ker(z)$ or $\ker(y+z)$.
But note that, since $\BBk =  \BBF_2$, we have
$$ \ker (z) = \{(0,0,0),(1,0,0),(0,1,0),(1,1,0)\}\subset \ker (x) \cup \ker (y) \cup \ker (x+y)$$
and 
$$ \ker (y+z) = \{(0,0,0),(1,0,0),(0,1,1),(1,1,1)\}\subset \ker (x) \cup \ker (y) \cup \ker (x+y).$$
This contradicts Definition \ref{def:mf}(b). So $\pi_3$ cannot be a part of an MAT-partition of $\CA$. Therefore, $\CA$ is not MAT-free after all, as claimed.
\end{example}

We see from Example \ref{ex:nonmf} that condition (b) in Definition \ref{def:mf} can fail for small fields. It is clear that this is not a lattice condition. If we replace it by a suitable lattice condition, we obtain a strictly larger class of arrangements capturing this notion of freeness independent of the underlying field.
Note that the two remaining conditions in Definition \ref{def:mf} (and in Theorem \ref{Thm_MAT}) are combinatorial. 

Specifically, if we replace condition (2) in Theorem \ref{Thm_MAT} by
\begin{itemize}
	\item [(2*)] $X \not \subseteq H$ for each $H \in \CA'$,
\end{itemize}
and condition (b) in Definition \ref{def:mf} by
\begin{itemize}
	\item [(b*)] $\cap_{H \in \pi_{k+1}} H   
	\nsubseteq H'$ for each $H' \in \Ac_k$,
\end{itemize}
then we obtain
a  larger class of arrangements, as obviously (2) implies (2*) and (b) implies (b*).  Denote this class by  $\CMF^*$. 
Both arrangements $\CA$ and $\CB$ from Example \ref{ex:nonmf} belong to $\CMF^*$. In particular, the example shows that 
$\CMF \subsetneq \CMF^*$.

\begin{defn}
	\label{def:mf*}
	The class $\CMF^*$ is the class of \emph{MAT${}^*$-free} arrangements.  
\end{defn}

Not only is condition (2*) easier to check than (2), more importantly, it only depends on the combinatorial data of the lattice of the arrangement ($X \wedge H \not=H$ for all $H \in \CA'$).
We thus readily derive the following.

\begin{theorem}
	\label{thm:main2}
The class $\CMF^*$ is combinatorial.
\end{theorem}

Remark \ref{rem:main}(i) gives Theorem \ref{thm:main}(ii).

\begin{remarks}
	\label{rem:main}
(i). If $\CA$ is defined over the infinite field $\BBk$, then $\CA$ is MAT-free if and only if it is MAT${}^*$-free. This follows from the prime avoidance theorem 
\cite[Lem.~3.3]{eisenbud}.

(ii). Furthermore, if $\CA$ is MAT${}^*$-free over a field extension of $\BBk$ then $\CA$ is also additively free over the smaller field $\BBk$, since additive freeness and MAT${}^*$-freeness are both combinatorial. 

(iii). Since we are predominantly interested in freeness as such, MAT${}^*$-freeness appears to be  a natural extension of of MAT-freeness.

(iv). There is a generalization of MAT-freeness, so called  \emph{MAT2-freeness}, see \cite[Def.~3.4]{CunMue19_MATfree}.
The discussion above about the (non-)combinatorial nature applies equally well to MAT2-freeness.
\end{remarks}

{\bf Acknowledgments}:
This work was supported by DFG-Grant
RO 1072/24-1 (DFG Project number 539858198) to G.~R\"ohrle.


\bigskip

\bibliographystyle{amsalpha}

\newcommand{\etalchar}[1]{$^{#1}$}
\providecommand{\bysame}{\leavevmode\hbox to3em{\hrulefill}\thinspace}


\end{document}